\documentclass[11pt]{article}
\usepackage{url} 
\usepackage{amsmath, amssymb}
\usepackage{amsthm}
\usepackage{setspace}
\usepackage[margin=1in]{geometry}
\usepackage[utf8]{inputenc}
\usepackage{algorithm}
\usepackage[noend]{algpseudocode}
\usepackage{graphicx}
\graphicspath{ {images/} }

\begin{document}

\title{On SAT Solvers and Ramsey-type Numbers}
\author{ 
Burcu Canakci\thanks{Microsoft Research at Cambridge, UK}\\
\small Microsoft Research \\[0.5ex]
\small \texttt{burcanakci@gmail.com}
\and
Hannah Christenson\thanks{Asana}\\
\small{Asana}\\[0.5ex]
\small \texttt{ } 
\and 
Robert Fleischman\thanks{TASH}\\
\small{TASH}\\[0.5ex]
\small \texttt{robbyef98@gmail.com} 
\and
William Gasarch\thanks{University of MD at College Park, Dept of CS}\\
\small{Univ of MD}\\[0.5ex] 
\small \texttt{gasarch@umd.edu} 
\and 
Nicole McNabb\thanks{Walmart Global Tech}\\
\small{Wallmart Global Tech}\\[0.5ex]
\small \texttt{ } 
\and 
Daniel Smolyak\thanks{University of MD at College Park, Dept of CS}
\small{Univ of MD}\\[0.5ex] 
\small \texttt{dsmolyak@gmail.com} 
}

\date{}
\maketitle

\begin{abstract}
We created and parallelized two SAT solvers to find new bounds on some Ramsey-type numbers. For $c > 0$, let $R_c(L)$ be the least $n$ such that for all $c$-colorings of the $[n]\times [n]$ lattice grid there will exist a monochromatic right isosceles triangle forming an $L$. Using a known proof that $R_c(L)$ exists we obtained $R_3(L) \leq 2593$. We formulate the $R_c(L)$ problem as finding a satisfying assignment of a boolean formula. Our parallelized probabilistic SAT solver run on eight cores found a 3-coloring of $20\times 20$ with no monochromatic $L$, giving the new lower bound $R_3(L) \geq 21$. We also searched for new computational bounds on two polynomial van der Waerden numbers, the "van der Square" number $R_c(VS)$ and the "van der Cube" number $R_c(VC)$. $R_c(VS)$ is the least positive integer $n$ such that for some $c > 0$, for all $c$-colorings of $[n]$ there exist two integers of the same color that are a square apart. $R_c(VC)$ is defined analogously with cubes. For $c \leq 3$, $R_c(VS)$ was previously known. Our parallelized deterministic SAT solver found $R_4(VS)$ = 58. Our parallelized probabilistic SAT solver found $R_5(VS) > 180$, $R_6(VS) > 333$, and $R_3(VC) > 521$. All of these results are new.
\end{abstract}

\section{Introduction}
\indent The following is often the first theorem taught in a lecture on Ramsey Theory. If there are six people at a party and for each pair of people either they know each other or neither knows the other, then there are either three people who all know each other or three people who all do not know each other. The reader can prove this by looking at cases.
\newline \indent Is this condition also true for a five-person party? How many people must we have to 
ensure that either four people all know each other or all don't know one another? Let $n$ represent 
the total number of people at the party. Let $K_n$ be the graph on $n$ vertices where every pair of 
vertices has an edge between them. Note that the party problem stated above can be restated as follows: 
for every 2-coloring of the edges of $K_6$ there exists a set of three points such that all of the edges 
between them are the same color.\\
\newline $\textbf{Notation.} $ $ \mathbb{N}$ is the set of natural numbers. If $n \in \mathbb{N}$, let $[n]$ be the set $\{1, 2, . . . n\}$. \\
\newline $\textbf{Definition.} $  Let $K_n$ be a complete graph on $n$ vertices. Let $V = [n]$ be the vertices of $K_n$ and $E$ the edges of $K_n$. Let $c \in [n]$. A $c$-coloring of the edges of $K_n$ is a function $f: E\rightarrow [c]$. \\
\newline $\textbf{Definition.} $  Let $G$ be a subgraph of $K_n$. If in a coloring of the edges of $K_n$ there is a monochromatic copy of $G$, then we will just say that there is a monochromatic $G$. We can again restate out party problem as follows: for every 2-coloring of the edges of $K_6$ there is a monochromatic $K_3$. \\

This party problem introduces the premise of Ramsey's Theorem,
initially proven by Ramsey~\cite{Ramsey-1930} (though easier proofs are
easily found on the web) 
which we state now:

\noindent
$\textbf{Theorem 1.1.}$ Ramsey's Theorem. 

Let $a$ and $b$ be positive integers. Then there exists a smallest positive integer $R(a,b)$ such that all 2-colorings of a complete graph on $R(a,b)$ vertices contain a monochromatic red $K_{a}$ or a monochromatic blue $K_{b}$. \\

We will refer to $R(a,b)$ as the Ramsey number for $(a,b)$. If $a = b$, we will just refer to the Ramsey number as $R(a)$. \\

The general theme of Ramsey Theory is that if you finitely color a big enough object
then some monochromatic subobject of a given size will always happen. 

Van der Waerden's theorem is in this spirit.
It was initially proven by 
Van der Waerden~\cite{VDW-1927} (though proofs in English are easily found on the web). 
and are closely related to the numbers we will discuss. 

\noindent
$\textbf{Theorem 1.2.}$ Van der Waerden's Theorem. \\
Let $k$ and $c$ be fixed positive integers. Then there exists a smallest positive integer $n$ such that all $c$-colorings of $[n]$ contain a monochromatic arithmetic sequence of integers of the form $\{a, a+d, a+2d, . . . a+(k-1)d\}$ for some positive integers $a$ and $d$ less than $n$. \\
\newline We will refer to $W(c,k)$ as the van der Waerden number with respect to $c$ and $k$. \\
\newline \indent In this paper, we will discuss two variants of Ramsey numbers we studied. We will introduce previously published theorems and results for each Ramsey-type number. Following previous results, we will present our new computational and theoretical results and current bounds.

\section{Theoretical Bounds on the $L$ Number}
$\textbf{Definition. }$ Let $n \in \mathbb{N}$. Consider an $[n] \times [n]$ grid of lattice points. An $L$ is a subgrid of the form $L = {(i, j), (i, j + t), (i + t, j + t)}$, or the isoscoles right triangle formed by the top left, bottom left, and bottom right corners of any square in the grid. \\
\newline $\textbf{Theorem 2.1.}$ $L$ Number Theorem. \\
Let $c$ be a positive integer. Then there exists a positive integer $n$ with $n > c$ such that all $c$-colorings of an $[n] \times [n]$ grid of lattice points contain a monochromatic $L$. \\
\newline We will refer to the integer $n$ guaranteed by the theorem above as the $L$ number, and denote it as $R_{c}(L)$ for a given $c$. \\
\newline \indent The first proof of the $L$ number theorem was given independently 
by Gallai (as reported by Rado~\cite{Rado-1933}, \cite{Rado-1943}) and Witt \cite{Witt-1951}
as a corollary to a more general theorem on Ramsey numbers. This proof gave an 
enormous upper bound for $R_{c}(L)$. A lower upper bound was found by 
Graham and Solymosi~\cite{GS-2006}.
The proof for Theorem 2.3, given below, was obtained through careful analysis 
of Graham and Solymosi's proof. Alexnovich \&  Manske~\cite{AM-2008} 
proved the specific case $R_2(L)$ 
through careful case analysis. There results is Theorem 2.2. 

\noindent
\textbf{Theorem 2.2.} 
$R_{2}(L) = 5$. 

\noindent
\textbf{Theorem 2.3.} 
$R_{3}(L)\leq 2593$. 

\noindent
{\bf Proof}  For any $3$-colored $[n]\times[n]$ grid, by the Pigeon-hole Principle there exists a color that occurs at least $\frac{n}{3}$ times along the diagonal. Without loss of generality, let this color be red. If there is no monochromatic $L$, then there must be ${\frac{n}{3} \choose 2} = \frac{n(n-3)}{18}$ points below the diagonal that are not colored red so as to avoid conflict with the red points on the diagonal. Let these points be colored blue and green. \\
\newline Since there are $n-1$ subdiagonals below the main diagonal, there exists a subdiagonal that contains at least $\frac{n(n-3)}{18(n-1)}$ points that are one of $c_2$ or $c_3$. Therefore, there are at least $\frac{n(n-3)}{36(n-1)}$ points on this diagonal that are all blue or all green.  Without loss of generality, assume that there are $\frac{n(n-3)}{36(n-1)}$ points colored blue. In order to continue to avoid a monochromatic $L$, there are $\frac{n(n-3)}{36(n-1)} \choose 2$ points below this diagonal that must be green, as they cannot be red or blue by our coloring of higher diagonals. \\
\newline Therefore, there exists a subdiagonal with $\frac{{\frac{n(n-3)}{36(n-1)} \choose 2}}{n-2}$ points that must be green, since there are at most $n-2$ subdiagonals below the second diagonal in consideration. However, if this number is greater than $1$, then there is a subdiagonal with at least $2$ points colored green, in which case the corner of the $L$ that they help form cannot be red, blue or green to avoid a monochromatic $L$. Because that point must be colored in some way, if $\frac{{\frac{n(n-3)}{36(n-1)} \choose 2}}{n-2} > 1$, then it is guaranteed that there is a monochromatic $L$ in the grid. If we take into account that we can round up to the nearest whole number when working with a number of points, the smallest $n$ such that this is the case is $2593$.

\section{New Computational Bounds on $L$}
\indent To our knowledge, no computational lower bounds on $R_{3}(L)$ were previously known. In order to find lower bounds on the $3$-colored $L$ number, we created two $k$-SAT solvers. Both solvers create a boolean formula 
in 3-CNF form by indexing the lattice points of an $[n]\times[n]$ grid. The first solver evaluates the formula and attempts to create a satisfying assignment deterministically using the framework of the classic DPLL algorithm. The basic steps of our algorithm are outlined in the pseudocode below. When any variable is set, the formula is reduced accordingly. 2SAT is an asymptotically linear time 2-SAT algorithm. Technically, the "most common variable" that the algorithm chooses is both weighted and considers positive and negative occurences of a variable separately. This means that we keep track of a separate "weight" for both $x$ and $\neg x$. Each time $x$ occurs, if the size of the clause we are considering is $s$, we add $\frac{1}{s}$ to the weight of $x$.  We then choose the literal with the highest weight as the "most common variable." This yielded the best results in our comparisons of a few different heuristics. 
\begin{algorithm}[H]
\caption{DPLL(formula f, partial assignment a)}
\begin{algorithmic}[1]
\For{$ \textit{var} \textbf{ in } \text{unit clause} \textbf{ in } \text{f}$} 
 \State $\textit{a}\text{[}\textit{var}\text{]} \gets true$
 \If {$\textit{clause} \textbf{ in }\textit{f}\text{ becomes unsatisfiable}$}
 \Return false
 \EndIf
\EndFor
\If {$\textit{f}\text{ is satisfied}$} \Return true
\EndIf
\If {$\textit{f} \textbf{ in } \text{2-CNF}$} \Return $\text{2SAT(f,a)}$ \EndIf
\For{$\textit{clause} \textbf{ in } \textit{f}$}
 \State $\textit{size} \gets \text{length of } \textit{clause}$
 \If {$\textit{size} = \text{2}$}
  \State $\textit{a}\text{[}\textit{clause}\text{[0]]} \gets true$
  \If {$\text{DPLL(f, a)}$}
  \State \Return true
  \Else \State $\textit{a}\text{[}\textit{clause}\text{[0]]} \gets      false$
  \State \Return $\text{DPLL(f,a)}$\EndIf
 \EndIf
\EndFor
$\textit{var} \gets \text{most common variable} \textbf{ in } \textit{f} $
\State $\textit{a}\text{[}\textit{var}\text{]} \gets true$
  \If {$\text{DPLL(f, a)}$}
  \State \Return true
  \Else \State $\textit{a}\text{[}\textit{var}\text{]} \gets false$
  \State \Return $\text{DPLL(f,a)}$\EndIf
\end{algorithmic}
\end{algorithm}

By running this serial recursive algorithm on a single processor, dual-core computer, we found a 3-coloring of a $14 \times 14$ lattice grid without a monochromatic $L$. Thus, our computational lower bound given by our serial DPLL algorithm is $n \geq 15$. \\ 
\indent Our second $k$-SAT solver uses a local search algorithm 
based on a description of the {\bf Walksat} algorithm. 
The algorithm generates a random assignment, and then flips a predetermined number of variables in an attempt to find a satisfying assignment. This algorithm is run a calculated large number of times in the solver. If an assignment is not found, the solver indicates that the formula is unsatisfiable. Note that since local search is a probabilistic algorithm, even if a satisfying assignment exists, the solver may not find it and will return unsatisfiable. The pseudocode below outlines the local search algorithm.
\begin{algorithm}[H]
\caption{LocalSearch(formula f, int maxflips, double p)}
\begin{algorithmic}[1]
\State $\textit{a} \gets \text{random coloring assignment}$
\State $\textit{u} \gets \text{initial unsatisfied clauses of }\textit{f}\text{ evaluated on }\textit{a}$
\For{$\textit{maxflips}\text{ times}$}
\If {$\textit{u} = \text{empty}$} \Return true \EndIf
\State $\textit{x} \gets \text{rand(0,1)}$
 \If {$\textit{x} < \textit{p}$}
  \State $\textit{var} \gets \text{random variable from an unsatisfied clause}$
 \Else
  \State $\textit{var} \gets \text{variable with "best flipping effect"}$
 \EndIf
 \State $\textit{a}\text{[var] = }\neg\textit{a} \text{[var]}$
 \State  $\textit{u} \gets \textit{u} \text{ updated with recent flip}$
\EndFor
\State \Return \textit{u}\text{ = empty}
\end{algorithmic}
\end{algorithm}

\indent The variable with the "best flipping effect" is the variable that, when flipped, results in the greatest decrease in the number of unsatisfied clauses. The local search algorithm as originally implemented performed slightly better than our implementation of DPLL. However, many clauses of the $L$ formula only serve to ensure that exactly one color could be assigned to each lattice point in the grid. Each of the $n^{2}$ points in the grid corresponds to $3$ variables: one representing red, one representing blue, and one representing green. Exactly one of these variables must be true in a satisfying assignment that represents a valid coloring of the grid. In light of this, we modified our algorithm in a few ways. \\
\indent First, we modified each randomly generated assignment so that only the clauses relating to the presence of a monochromatic $L$, and not those that guarantee a valid coloring, could be unsatisfied. In addition, when we flip a variable, we flip one of the corresponding variables for that spot in the grid in order to maintain the invariant that the assignment respresents a valid coloring. We take into account the net effect of these multiple flips when choosing the variable, and when we have an option of which color-representing variable to flip, we also pick the "better" one by this metric. This modified local search algorithm, run on the same single processor computer, found a satisfying assignment for $n = 19$, though it returned unsatisfiable for all of the tests that we ran on $n = 20$. Because local search is a probabilistic algorithm, its inability to find a satisfiying assignment for higher values of $n$ does not imply that they do not exist. \\
\indent Next, we parallelized both solvers because we were able to run them on an 8-core virtual machine cluster. We parallelized the local search algorithm by simply running local search many times at once, with different random seeds, allowing the solver to search for assignments in parallel. In the deterministic solver, we created explicit threads at the first few steps of the DPLL recursion, which then completed in parallel. Our parallelized DPLL algorithm only found a satisfying assignment for $n = 15$, which was slightly better than our result from our serial DPLL algorithm, though not as good as the result achieved by our serial local search algorithm. Our parallelized local search algorithm found satisfying assignments for $n = 20$, one of which is shown in Figure 1 below. Therefore, our current lower bound for the $L$-number is $21$.
\newline 
\newline $\textbf{Theorem 3.1.}$ \\
$21 \leq R_{3}(L) \leq 2593$. \\
\newline $Proof.$  By Theorem 2.3, $R_{3}(L) \leq 2593$. Our parallelized local search algorithm run on an 8-core virtual machine cluster found satisfiable 3-colored $20 \times 20$ lattice grids. Thus, $R_{3}(L) \geq 21$ is our lower bound by computation.
\newline
\newline {\bf Figure 1: A 20 $\times$ 20 grid with no monochromatic $L$.} 
\newline \includegraphics{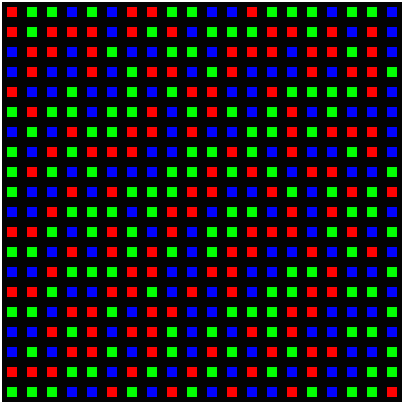}

\section{Theoretical Bounds on Polynomial van der Waerden Numbers}
$\textbf{Theorem 4.1.}$ General Polynomial van der Waerden Theorem. \\
Let $c \in \mathbb{N}$. For all sets of polynomials $\{p_{1}, . . . p_{k}\}$ in $\mathbb{Z}[x]$ with a zero constant term, there exists an $n \in \mathbb{N}$ where $n > c$ such that for all $c$-colorings of $[n]$, there exist positive integers $a$ and $d$ that guarantee a monochromatic arithmetic sequence of the form $\{a, a+p_{1}(d), . . . a+p_{k}(d)\}$. \\
\newline \indent The above theorem, which generalizes van der Waerden's theorem to polynomials with zero constant term, was first proven by Bergelson and Liebman~\cite{BL-1999}.
Walters~\cite{Walters-2000} gave the first constructive proof 
for the polynomial van der Waerden theorem, 
which gave explicit, albeit extremely large, upper bounds. 
In this paper, we focus on two specific polynomial van der Waerden numbers introduced in the corollaries below. \\
\newline
$\textbf{Corollary 4.1.}$ Van der Square Numbers. \\
Let $c \in \mathbb{N}$. Then there exists an $n \in \mathbb{N}$, $n > c$, such that every $c$-coloring of $[n]$ contains two monochromatic integers $x^{2}$ distance apart for some integer $x$. \\
\newline
$\textbf{Corollary 4.2.}$ Van der Cube Numbers. \\
Let $c \in \mathbb{N}$. Then there exists an $n \in \mathbb{N}$, $n > c$, such that every $c$-coloring of $[n]$ contains two monochromatic integers $y^{3}$ distance apart for some integer $y$. \\
\newline \indent The proofs for these corollaries follow directly from Theorem 4.1. We will denote the van der Square number defined in Corollary 4.1 by $R_{c}(VS)$ for a given $c$. Similarly, we will denote the van der Cube number defined in Corollary 4.2 by $R_{c}(VC)$ for a given $c$. \\
\newline 
$\textbf{Theorem 4.2.}$ \\
$R_{2}(VS) = 5$. \\
\newline $Proof.$ First, let us prove $R_{2}(VS) \leq 5$. Without loss of generality, assume $1$ is colored red. Then each subsequent number in the sequence must be colored the opposite of the number before it, as consecutive numbers are $1^{2}$ apart. Thus, we color $2$ blue, $3$ red, $4$ blue, and $5$ red. However, $5$ is $2^{2}$ apart from $1$ and they are both colored red. If $5$ was blue, $5$ and $4$ would be the same color, and they are separated by $1^2$. Therefore, in any $2$-coloring of $[5]$, there is a monochromatic pair of integers separated by a square, so $R_{2}(VS) \leq 5$.
\newline \indent The table below shows a possible $2$-coloring of a length $4$ sequence with no monochromatic pair of numbers a square apart. Thus, $R_{2}(VS) = 5$.
\begin{table}[h!]
  \begin{center}
    \caption{Sample 2-coloring of a sequence of length 4}
    \begin{tabular}{c|c|c|c}
      1 & 2 & 3 & 4 \\
      \hline
      R & B & R & B \\
    \end{tabular}
  \end{center}
\end{table} \\
\newline $\textbf{Theorem 4.3.}$ \\
$R_{3}(VS) \leq 68$. \\
\newline $Proof.$
We will attempt to color a sequence of $1$ through $68$ without a monochromatic pair of numbers separated by a square, in order to show that such a task is impossible. Without loss of generality, assume $10$ is colored red. Then $1$ and $26$ cannot be red because they are a square apart from $10$. $1$ and $26$ also cannot be the same color since they are a square apart, thus, without loss of generality, color $1$ blue and color $26$ green. Next, $17$ is a square apart from both $1$ and $26$, forcing it to be red. When this process is repeated, all numbers that are $7$ apart from $10$ to $59$, $\{10, 17, 24, . . ., 59\}$ must be red. However, $10$ and $59$ are $7^2$ apart. $59$ was forced to be red by the forced green and blue color of the numbers $68$ and $43$. Therefore, $R_{3}(VS) \leq 68$.
\begin{table}[h!]
  \begin{center}
    \caption{Forced number colorings in any 3-coloring}
    \begin{tabular}{c|c|c|c}
      1 & 10 & 17 & 26\\
      \hline
      B & R & R & G\\
    \end{tabular}
  \end{center}
\end{table}

\section{New Computational Bounds on Polynomial van der Waerden Numbers}
\indent Using our $k$-SAT solvers with the heuristics discussed in a previous section, we found computational results and lower bounds for some polynomial van der Waerden numbers. The table below summarizes the current known values and bounds on $R_{c}(VS)$. The bolded values represent the results given by our $k$-SAT solvers, which were previously unknown. $R_{c}(VS)$ for $c \leq 3$ were previously known.
\begin{table}[h!]
  \begin{center}
    \caption{Known van der Waerden square numbers}
    \begin{tabular}{|c||c|c|c|c|c|c|}
       c & 1 & 2 & 3 & 4 & 5 & 6 \\
      \hline
      n & 2 & 5 & 29 & \textbf{58} & $>$\textbf{180} & $>$\textbf{333} \\
    \end{tabular}
  \end{center}
\end{table}

In addition, our parallelized local search algorithm was able to find that $R_{3}(VC) > 521$ using the heuristics discussed earlier in the paper. It is easy to see that $R_2(VC) = 9$ in a proof parallel to that presented for $R_2(VS) = 4$. However, our DPLL algorithm was not able to provide reasonable bounds for $R_{3}(VC)$ or other van der Cube numbers for higher values of $c$. In the future, we hope to continue to improve both of our parallelized algorithms and to gain additional computational power in order to compute tighter bounds on polynomial van der Waerden numbers in addition to still unknown Ramsey and $L$ numbers.


\end{document}